\documentclass{rspublic}

\begin{document}

\title[Commutativity up to a Factor]{Commutativity up to a Factor of Bounded Operators in Complex Hilbert Space}

\author[J.A. Brooke, P. Busch, D.B. Pearson]{James A. Brooke$^{(1)}$, Paul
Busch$^{(2)}$, and David B. Pearson$^{(2)}$}

 \affiliation{
$^{(1)}$Department of
Mathematics and Statistics, University of Saskatchewan, Saskatoon, SK S7N 5E6, Canada\\
{\small E-mail: brooke@usask.ca}\\
$^{(2)}$Department of Mathematics, University of Hull, Hull HU6
7RX, UK\\
{\small E-mail: p.busch@maths.hull.ac.uk,
d.b.pearson@maths.hull.ac.uk}}

\label{firstpage}

\maketitle

\begin{abstract}{Hilbert space operators, commutator, commutation
relations, self-adjointness, positivity, quantum measurement
theory} We explore commutativity up to a factor, $AB=\lambda BA$,
for bounded operators in a complex Hilbert space. Conditions on
the possible values of the factor $\lambda$ are formulated and
shown to depend on spectral properties of the operators involved.
Commutativity up to a unitary factor is considered for pairs of
self-adjoint operators. Examples of nontrivial realizations of
such commutation relations are given.
\end{abstract}

\section{Introduction}

Commutation relations between self-adjoint operators in a complex
Hilbert space $H$ play an important role in the interpretation of
quantum mechanical observables and the analysis of their spectra.
Accordingly such relations have been extensively studied in the
mathematical literature (see, e.g., the classic study of Putnam
(1967)). An interesting related aspect concerns the commutativity
up to a factor for pairs of operators. Certain forms of
noncommutativity can be conveniently phrased in this way. This is
the case with the famous canonical (or Heisenberg) commutation
relations for position $Q$ and momentum $P$,
\begin{equation*}
QP-PQ\subset \ri I,
\end{equation*}
which can be recast in the form of the Weyl relations,
\begin{equation*}
\exp (\ri\alpha Q)\exp (\ri\beta P)=\exp (-\ri\alpha \beta )\exp
(\ri\beta P)\exp (\ri\alpha Q)\;\;[\alpha ,\beta \in \mathbb{R}].
\end{equation*}
Another example well known in the physical context are
anticommutation relations between (non-selfadjoint) fermionic
creation and annihilation operators, or between Pauli spin
matrices in $\mathbb{C}^{2}$, e.g., $\sigma _{x}\sigma
_{y}=\ri\sigma_z=-\sigma _{y}\sigma _{x}$, where
\begin{equation*}
\sigma _{x}=\left(
\begin{array}{ll}
\text{0} & \text{1} \\
\text{1} & \text{0}
\end{array}
\right) ,\;\;\sigma _{y}=\left(
\begin{array}{ll}
{0} & {-}{\ri} \\
{\ri} & {0}
\end{array}
\right) ,\;\;\sigma_z=\left(
\begin{array}{ll}
{1} & {0} \\
{0} & {-1}
\end{array}
\right).
\end{equation*}

More recently, algebraic relations of the form $yx=q\,xy$ have
been the subject of study in the context of quantum groups (e.g.,
Kassel 1995), and their matrix realizations give examples of
operator pairs commuting up to a
factor. The \emph{quantum enveloping algebra }$U_{q}\left( \frak{sl\,}%
(2)\right) $ of the $\frak{sl\,}(2)$ Lie algebra is defined, for $q$ with $%
q^{2}\neq 1$, to be the algebra with four generators $E,F,K,K^{-1}$ subject
to the relations
\begin{eqnarray*}
KK^{-1} &=&K^{-1}K=I, \\
KEK^{-1}=q^{2}E,\; &&KFK^{-1}=q^{-2}F, \\
EF-FE &=&\frac{K-K^{-1}}{q-q^{-1}}.
\end{eqnarray*}
An explicit realization is given later.

In this paper we consider pairs of bounded operators $A,B$ on a complex
Hilbert space $H$ and explore the conditions under which they can commute up
to a factor, i.e.,
\begin{equation*}
AB=\lambda BA,\;\;\lambda \in \mathbb{C}\setminus\{0\}.
\end{equation*}
Obviously, if $AB=O$ then $BA=O$, that is the operators commute
and $\lambda $ can be any (nonzero) complex number. The remaining
case of $AB\neq O$ is dealt with in a sequence of special cases.
We also present a result concerning the more general problem of
commutativity up to a unitary factor.

The main results are the following.

\begin{theorem}
\label{theo1}Let $A,B$ be bounded operators such that $AB\neq O$ and $%
AB=\lambda BA$, $\lambda \in \mathbb{C}$. Then:\newline (i) if $A$
or $B$ is self-adjoint then $\lambda \in \mathbb{R}$;\newline (ii)
if both $A$ and $B$ are self-adjoint then $\lambda \in \left\{
1,-1\right\} $;\newline (iii) if $A$ and $B$ are self-adjoint and
one of them positive then $\lambda =1$.
\end{theorem}

\begin{theorem}
\label{theo3}Let $A,B$ be bounded and self-adjoint. The following
are equivalent:\newline (i)$\;\ \ AB=UBA$ for some unitary
$U$.\newline (ii) $AB^{2}=B^{2}A$ and $BA^{2}=A^{2}B$.\newline If
$A$ or $B$ is positive then (i) is equivalent to $AB=BA$.
\end{theorem}

Section 2 provides a proof of theorem \ref{theo1}  together with
some further case studies highlighting the interrelations between
the spectral properties and admissible values of the factor
$\lambda$. Section 3 offers a substantially different proof
technique for part (iii) of theorem \ref{theo1}. The case of
commutativity up to a unitary factor will be discussed in \S4 for
pairs of self-adjoint operators, leading to a proof of theorem
\ref{theo3}. The paper concludes with a number of realizations of
commutation relations up to a factor (\S5) and a sketch of a
quantum mechanical application (\S6).

\section{Spectral properties and conditions on $\protect\lambda $}

We will use repeatedly the following well known result concerning
the spectrum of a product of two bounded operators $X,Y$:
\begin{equation*}
\sigma (XY)\backslash \{0\}=\sigma (YX)\backslash \{0\}
\end{equation*}
(e.g., Halmos 1982, Problem 76). From this we obtain the following
useful observations.

\begin{lemma}
\label{pro0}Let $AB=\lambda BA$ and $AB\neq O$. Then $0$ is either
in both or neither of $\sigma \left( AB\right) $ and $\sigma
\left( BA\right) $. Hence
\begin{equation*}
\sigma \left( AB\right) =\sigma \left( BA\right) =\lambda \sigma \left(
AB\right) .
\end{equation*}
If $0\in \sigma \left( AB\right) $, at least one of $A$ or $B$ does not have
a bounded inverse. If $\sigma \left( AB\right) \neq \left\{ 0\right\} $,
then $\left| \lambda \right| =1$. If $0\notin \sigma \left( AB\right) $,
both $A$ and $B$ have bounded inverses.
\end{lemma}

\begin{proof} Observe that $\lambda \neq 0$ since $AB\neq O$. Suppose
$0\in \sigma \left( AB\right) $ and $0\notin \sigma \left(
BA\right) $. Then $BA$ is invertible, and $AB\left( BA\right)
^{-1}=\lambda I$. The operator on the left hand side has $0$ in
its spectrum and so $\lambda =0$.
This is a contradiction. Similarly, $0\notin \sigma \left( AB\right) $ and $%
0\in \sigma \left( BA\right) $ cannot occur. Hence, either $0\in \sigma
\left( AB\right) \cap \sigma \left( BA\right) $ or $0\notin \sigma \left(
AB\right) \cup \sigma \left( BA\right) $. It follows that $\sigma \left(
AB\right) =\sigma \left( BA\right) $.

If $0\in \sigma \left( AB\right) =\sigma \left( BA\right) $ then $0$ is in
the spectrum of $A$ or $B$. For otherwise both $A$ and $B$, and therefore $%
AB $, have bounded inverses, so that $0\neq \sigma \left( AB\right) $.

Let $\sigma \left( AB\right) \neq \left\{ 0\right\} $. The set $\sigma
(AB)=\lambda \sigma (AB)$ is compact and hence there is an element, $\gamma $%
, with maximal modulus. If $\left| \lambda \right| >1$ then $\lambda \gamma $
has modulus greater than $\left| \gamma \right| $ and hence does not belong
to $\sigma (AB)$. This shows that $\left| \lambda \right| \leq 1$. A similar
argument for $\lambda ^{-1}\sigma \left( AB\right) =\sigma \left( AB\right) $
shows that $\left| \lambda ^{-1}\right| \leq 1$, and so $\left| \lambda
\right| =1$.

If $0\notin \sigma \left( AB\right) =\sigma (BA)$, then $AB$ and $BA$ have
bounded inverses since their spectra are bounded away from $0$. Hence both $%
A $ and $B$ are surjective and injective and thus they have bounded
inverses.
\end{proof}

This result entails, in particular, that $AB=\lambda BA(\neq O)$, with $%
\left| \lambda \right| \neq 1$, can only be realized if $AB$ is
quasi-nilpotent. We are now ready to prove theorem \ref{theo1}.

\noindent \emph{Proof of Theorem \ref{theo1}.} (i): Let $A=A^{\ast
}$. The condition $AB=\lambda BA(\neq O)$ implies
\begin{equation*}
ABB^{\ast }=\lambda BAB^{\ast }=\lambda \overline{\lambda }^{-1}BB^{\ast }A.
\end{equation*}
This is to say that $ABB^{\ast }$ is a multiple of a self-adjoint operator,
with spectrum
\begin{equation*}
\sigma (ABB^{\ast })=\lambda \overline{\lambda }^{-1}\sigma (BB^{\ast
}A)=\lambda \sigma (BAB^{\ast })\subseteq \lambda \cdot \lbrack -\left\|
BAB^{\ast }\right\| ,\left\| BAB^{\ast }\right\| ].
\end{equation*}
By lemma \ref{pro0} we have $\sigma (ABB^{\ast })=\sigma (BB^{\ast }A)$%
. Furthermore, since $A$ and $BB^{\ast }$ are self-adjoint, $(ABB^{\ast
})^{\ast }=BB^{\ast }A$, and so
\begin{equation*}
\sigma (ABB^{\ast })=\sigma (BB^{\ast }A)=\overline{\sigma (ABB^{\ast })}%
\subseteq \lambda \cdot \lbrack -\left\| BAB^{\ast }\right\| ,\left\|
BAB^{\ast }\right\| ].
\end{equation*}
Multiplication by $\lambda \in \mathbb{C}$ scales the line segment $%
[-\left\| BAB^{\ast }\right\| ,\left\| BAB^{\ast }\right\| ]$ by a factor $%
\left| \lambda \right| $ and rotates it in the complex plane about
the origin; and the conjugation invariance of $\sigma (ABB^{\ast
})$ means symmetry of the set under reflection in the real axis.
Since $\sigma (BB^{\ast }A)=\lambda \sigma (BAB^{\ast })\neq
\left\{ 0\right\} $, it follows that $\lambda \in \mathbb{R}$ or
$\lambda \in \ri\mathbb{R}$.

The case $\lambda \in i\mathbb{R}$ leads to $\lambda \overline{\lambda }%
^{-1}=-1$, and so from the above, $ABB^{\ast }=-BB^{\ast }A$. Multiply both
sides with $A$ from the right to get $ABB^{\ast }A=-BB^{\ast }A^{2}$. The
left-hand side is non-negative; therefore the right-hand side must be
self-adjoint, hence $BB^{\ast }$ and $A^{2}$ commute, so their product is
both non-negative and non-positive. Thus, $BB^{\ast }A^{2}=O$, in
contradiction to the assumption $AB\neq O$. Hence $\lambda \not\in \ri\mathbb{R%
}$ and so $\lambda \in \mathbb{R}$.

(ii): We have $A=A^{\ast }$ and $B=B^{\ast }$. Then
\begin{equation*}
AB^{2}=\lambda BAB=\lambda ^{2}B^{2}A.
\end{equation*}
From (i) we have $\lambda \in \mathbb{R}$, so $\lambda BAB$, and hence $%
AB^{2} $, is self-adjoint, that is, $A$ commutes with $B^{2}$. Therefore $%
\lambda ^{2}=1$.

(iii): This follows directly from the last equation: if $A$ is positive then $%
AB^{2}$ and $BAB$ are positive as well, and so by (ii), $\lambda
=1$. \qed


We note that the condition of positivity is necessary: the Pauli
spin operators $\sigma _{x},\sigma _{y}$ acting in
$\mathbb{C}^{2}$, satisfy $\sigma _{x}\sigma _{y}=-\sigma
_{y}\sigma _{x}=\ri\sigma_z\neq 0$, but, as $\sigma _{y}^{2}=I$,
then $\sigma _{x}\sigma _{y}^{2}=\sigma _{x}=\sigma _{y}^{2}\sigma
_{x}$.

Another  interesting class of  cases arises if one of $A$ and $B$
has a bounded inverse.

\begin{proposition}
\label{prop-inv}Let $AB=\lambda BA$ and $AB\neq O$, and assume
that $A$ has a bounded inverse. Then
\begin{equation*}
\sigma (B)=\lambda \cdot \sigma (B).
\end{equation*}
Furthermore, if $\sigma (B)\neq \left\{ 0\right\} $ or $A$ is unitary then $%
\left| \lambda \right| =1$.
\end{proposition}

\begin{proof} We have $ABA^{-1}=\lambda B$. Since a similarity
transformation leaves the spectrum unchanged, we immediately
obtain $\sigma (B)=\lambda \cdot \sigma (B)$. Now suppose $\sigma
(B)\neq \left\{ 0\right\} $. A similar argument as in lemma
\ref{pro0} shows that $\left|
\lambda \right| \leq 1$ and $\left| \lambda ^{-1}\right| \leq 1$, and so $%
\left| \lambda \right| =1$. Finally, if $A$ is unitary then
$\left\| B\right\| =\left\| ABA^{-1}\right\| =\left| \lambda
\right| \left\| B\right\| $, so $\left| \lambda \right| =1$.
\end{proof}

The condition $\sigma (B)=\lambda \cdot \sigma (B)$ for given
$\lambda $  limits the choice of operators $B$ that could enter
the relation $AB=\lambda BA$. Indeed, multiplication with $\lambda
=\re^{\ri\theta }$ corresponds to a rotation about the origin in
the complex plane. Thus the spectrum of $B$
must be invariant under rotation about $0$ by the angle $\theta $. If $%
\lambda $ is not a primitive root of unity, then for any $\beta
\in \sigma \left( B\right) $, the set $\left\{ \lambda ^{n}\beta
:n=1,2,\dots \right\} $ will be dense in the circle $\left\{ z\in
\mathbb{C}:\left| z\right| =\left| \beta \right| \right\} $; then
the spectrum of $B$ must be the union of such circles. On the
other hand, if $\lambda ^{k}=1$ for some minimal $k\in
\mathbb{N}$, then $\sigma \left( B\right) $ is periodic with minimal period $%
\frac{2\pi }{k}$ under rotations about $0$.

If $B$ is given then the spectral condition $\sigma \left(
B\right) =\lambda \cdot \sigma \left( B\right) $ entails
constraints on the possible values of $\lambda $. For example, if
$\sigma (B)$ lies on a straight line segment originating at $0$,
then we must have $\lambda =1$. If it lies on a line segment which
contains $0$ in the interior, and if $\sigma (B)$ is not symmetric
under reflection at the origin then again $\lambda =1$; if $\sigma
(B)$ is symmetric under this reflection then $\lambda $ can be $1$
or $-1$.

In finite dimensional Hilbert spaces, trace or determinant
conditions give further simple constraints on the factor
$\lambda$.

\begin{proposition}
\label{prop-fin}Let $\dim H=n<\infty $, and assume $AB=\lambda BA$ and $%
AB\neq O$. If the trace $\mathrm{tr}\left[ AB^{k}\right] \neq 0$ or $\mathrm{%
tr}\left[ A^{k}B\right] \neq 0$ for some $k\in \mathbb{N}$, then $\lambda =1$%
. If $\det \left[ AB\right] \neq 0$ (i.e. $A,B$ are both invertible) then $%
\lambda ^{n}=1$.
\end{proposition}

\begin{proof} We have $AB^{k}=\lambda BAB^{k-1}$. On taking the
trace, this yields $\mathrm{tr}\left[ AB^{k}\right] =\lambda \mathrm{tr}%
\left[ AB^{k}\right] \neq 0$, hence $\lambda =1$. The second
statement follows similarly.
\end{proof}

\section{A proof based on resolvent integrals}

We provide an alternative proof of theorem \ref{theo1}, statement
(iii), based on resolvent techniques.

On taking adjoints in $AB=\lambda BA$, $BA=\overline{\lambda }AB$, which
leads to $\lambda \overline{\lambda }=1$. We consider first the case where $%
A $ has purely discrete spectrum. Let $f$ be an eigenvector of $A$ with
strictly positive eigenvalue $\alpha $. Then $Af=\alpha f$, so that $%
BAf=\alpha Bf=\lambda ^{-1}ABf$ and hence $A(Bf)=\alpha \lambda (Bf)$. If $%
\lambda =-1$ or $\Imag\lambda \neq 0$ then since $A$ cannot have a
negative or complex eigenvalue, it follows $Bf=0$; hence, in this case $%
BAf=\lambda ^{-1}ABf=0$ for any eigenvector of $A$ having strictly positive
eigenvalue. But we also have $BAf=0$ for any eigenvector of $A$ having zero
eigenvalue.

Since the eigenvectors of $A$ span the Hilbert space, we deduce
that $BA=O$ if $\lambda =-1$ or $\Imag\lambda \neq 0$. Since we
assumed $AB\neq O$, then we are left with $\lambda =1$.

General case: For any $z\in \mathbb{C}$ we have $(A-zI)B=\lambda B(A-\frac{z%
}{\lambda }I)$. For $z\notin \sigma (A)$ and $\frac{z}{\lambda }\notin
\sigma (A)$ we have $(A-zI)^{-1}B=\lambda ^{-1}B(A-\frac{z}{\lambda }I)^{-1}$%
. Let $\mathcal{J}=[a,b]$ be a closed bounded interval in $(0,\infty )$ and
consider the spectral projection $E^{A}(\mathcal{J})$ of $A$, associated
with $\mathcal{J}$. Assume the endpoints of $\mathcal{J}$ are not
eigenvalues of $A$. Then
\begin{equation*}
E^{A}(\mathcal{J})B=\frac{1}{2\pi \ri}\lim_{\epsilon \rightarrow
+0}\int_{a}^{b}\left[ \left( A-(t+\ri\epsilon )\right)
^{-1}-\left( A-(t-\ri\epsilon )\right) ^{-1}\right] \,\cdot
\,B\,dt
\end{equation*}
(where the weak limit is to be taken). So
\begin{equation*}
E^{A}(\mathcal{J})B=\frac{1}{2\pi \ri}\lim_{\epsilon \rightarrow
+0}\int_{a}^{b}\lambda ^{-1}B\,\,\cdot \left[ \left( A-\frac{t+\ri\epsilon }{%
\lambda }\right) ^{-1}-\left( A-\frac{t-\ri\epsilon }{\lambda }\right) ^{-1}%
\right] \,\,dt.
\end{equation*}
Assume $\lambda =-1$ or $\Imag\lambda \neq 0$. Then, for
$\epsilon$
small enough, the distance between $\frac{t\pm \ri\epsilon }{\lambda }$ and $%
\sigma (A)$, for $t\in \lbrack a,b]$, is bounded below by a positive number.
Hence
\begin{equation*}
\left\| \left( A-\frac{t\pm \ri\epsilon }{\lambda }\right)
^{-1}\right\| \leq \const,
\end{equation*}
and the integrand above is bounded in norm by $\const\,\frac{\mathrm{%
\epsilon }}{|\lambda |^{2}}$. We conclude that $E^{A}(\mathcal{J})B=O$ and
thus $E^{A}\left( [c,\infty )\right) B=O$ for any $c>0$. Any vector in the
range of $B$ can only be zero or an eigenvector of $A$ with eigenvalue zero.
\endproof

\section{Commutativity up to a unitary factor}

An obvious generalization of the preceding considerations concerns
commutativity up to a unitary factor. We provide a discussion in
the case of  $A,B$ being bounded and self-adjoint operators
(theorem \ref{theo3}).

\begin{proposition}
\label{uba}Let $A,B$ be bounded self-adjoint operators on $H$. The
following are equivalent:\newline (i) $\;AB^{2}A=BA^{2}B$;\newline
(ii)$\;AB=UBA$ for some unitary $U$.
\end{proposition}

\begin{proof}(ii)$\Longrightarrow $(i): Observe that $%
BA=ABU^{\ast }$, and so $AB=UABU^{\ast }$. Hence $AB$ commutes with $U$ and $%
U^{\ast }$, and similarly for $BA$. Thus we get
\begin{equation*}
AB^{2}A=UBAU^{\ast }AB=BAUU^{\ast }AB=BA^{2}B.
\end{equation*}
(i)$\Longrightarrow $(ii): The condition (i) is equivalent to
$|AB|=|BA|$, where we have $|C|=\left( C^{\ast }C\right) ^{1/2}$.
Thus $AB$ is normal. It follows that
\begin{equation*}
\ker AB=\ker BA=\ker |AB|=\ker |BA|.
\end{equation*}
Let $Q$ denote the associated projection. We then have that the
closures of the ranges of these four operators all coincide. Let
$P$ denote the associated projection. Hence $P+Q=I$. By the polar
decomposition theorem there are partial isometries $V,W$ such that
$AB=V|AB|$, $BA=W|BA|$. Note
that $VV^{\ast }=V^{\ast }V=WW^{\ast }=W^{\ast }W=P$. We have $%
AB^{2}A=V|AB|^{2}V^{\ast }$, so that $V,V^{\ast }$ commute with
$|AB|$. (This follows from $\left| AB\right| P=\left| AB\right|
$.) Now we have  $AB=V|AB|=|BA|W^{\ast }=W^{\ast }|AB|$, and since
the range of $\left| AB\right| $ is dense in $PH$, it follows that
$W^{\ast }=V$, and $V^{\ast }=W $. Finally,
\begin{equation*}
AB=V\left| AB\right| =VW^{\ast }BA=V^{2}BA.
\end{equation*}
We can extend $V^{2}|_{PH}$ to a unitary map $U:=V^{2}|_{PH}\oplus
Q$.
\end{proof}

Note that we have the following:
\begin{equation*}
AB=PAP\,PBP=APB.
\end{equation*}
Proof: For $f\in H$, consider $Bf=PBf+g$; then $g\in \ker AB$. It
follows
that $Ag\in \ker B\subseteq \ker AB=\ker BA$. Therefore, $ABf=APBf+Ag$, and $%
PAg=0$, so $PABf=PAPBf$. Replace $f$ with $Pf$ to get
\begin{equation*}
PABPf=ABf=PAPBPf.
\end{equation*}
This means that $AB=UBA$ is equivalent to $PAPPBP=UPBPPAP$, so
that it suffices to consider the case $Q=O$.

We will make use of the following result due to Gudder \& Nagy
(2000).

\begin{theorem}
\label{theo2}Let $A,B$ be bounded and self-adjoint. The following
are equivalent:\newline (i) \ \ $AB^{2}A=BA^{2}B$;\newline (ii)
$AB^{2}=B^{2}A$ and $BA^{2}=A^{2}B$.
\end{theorem}

\noindent The proof of this theorem is based on the
Fuglede-Putnam-Rosenblum theorem (Rudin 1973, Theorem 12.16,
Halmos 1982, Problem 152). Combination of proposition \ref{uba}
and theorem \ref{theo2} completes the proof of theorem
\ref{theo3}.\hfill$\square$

If $A$ or $B$ is positive it follows from $AB^2=B^2A$ and
$A^2B=BA^2$ that $A$ and $B$ commute. Hence in this case theorem
\ref{theo3} reduces to the case where $U$ can be chosen to be $I$,
in accordance with theorem \ref{theo1}. However, the following
example shows that theorem \ref{theo3} comprises a wider range of
cases than theorem \ref{theo1}; that is, $AB=UBA$ can be satisfied
for non-commuting $A$ and $B$ in such a way that $U$ cannot be one
of $I,-I$.

Let $H={\mathbb C}^2$, and use the Pauli spin operators
$\sigma_x,\sigma_y,\sigma_z$ to define the self-adjoint operators
$A=\sigma_x$, $B=\frac{1}{\surd 2}(\sigma_x+\sigma_y)$, and the
unitary operator $U=\ri\sigma_z$. Using the relations
$\sigma_x\sigma_y=-\sigma_y\sigma_x=\ri\sigma_z$,
$\sigma_x^2=\sigma_y^2=\sigma_z^2=I$, one readily verifies that
$A^2=B^2=I$, so that $AB^2=B^2A$ and $A^2B=BA^2$. Furthermore we
obtain $ \sigma_x(\sigma_x+\sigma_y)=I+\ri\sigma_z$,
$(\sigma_x+\sigma_y)\sigma_x=I-\ri\sigma_z$, and so
$\ri\sigma_z(I-\ri\sigma_z)=I+\ri\sigma_z$. Hence we see that
$AB=UBA$ with a non-trivial unitary $U$.


We note that boundedness is an essential requirement in theorem
\ref{theo3}: the operators of position $Q$ and momentum $P$ do
satisfy $QP^{2}Q=PQ^{2}P$ on a dense domain, but $QP^{2}\neq
P^{2}Q$. Thus theorem \ref{theo3} does not extend to unbounded
operators. However, we note that for pairs of (not necessarily
bounded) self-adjoint operators $A,B$ whose products $AB$, $BA$
exist on dense domains, the polar decomposition theorem for closed
operators can be applied. Therefore we conjecture that the
statement of proposition \ref{uba} can be adjusted such as to
cover these cases.

\section{Realizations}

We shall give some realizations of $AB=\lambda BA$.

\vspace{3pt}

\noindent{\bf Example 1.} An example for the case where both $A$
and $B$ are unitary is given by the Weyl commutation relations
mentioned in the introduction. The essence of this case can be
conveniently exhibited in the following example. Let $H=L^{2}(S^{1},d\phi )$%
, and define $Bf(\phi )=e^{i\phi }f(\phi )$, $\phi \in \left[ 0,2\pi \right]
$. Here $\sigma \left( B\right) =\left\{ \lambda \in \mathbb{C}:\left|
\lambda \right| =1\right\} $. This set is invariant under arbitrary
rotations about $0$, that is $\sigma (B)=\lambda \cdot \sigma (B)$ for any $%
\lambda =\re^{\ri\theta }$. A unitary operator $A$ that permutes
the spectral
projections of $B$ in a suitable way is given by $Af(\phi )=f(\phi +\theta )$%
. Then $ABf(\phi )=\re^{\ri\theta }BAf\left( \phi \right) $.


\vspace{3pt}

\noindent{\bf Example 2.} For pairs of unitary operators $A,B$ \
in a finite-dimensional Hilbert space, by proposition
\ref{prop-fin}, $\lambda ^{n}=1$ in accordance with the fact that
the spectra are finite; and a similar construction to the
preceding example can be carried out to realize $B$ as a
multiplication operator and $A$ as a cyclic shift. This case was
treated by Weyl (1928).


\vspace{3pt}


\noindent{\bf Example 3.} In Problem 238 of Halmos (1982), it is
shown that in an infinite-dimensional Hilbert space $H$,
$AB=\lambda BA$ can be realized for
any $\lambda $ of modulus 1 with a pair of unitary operators $A,B$. That $%
\left| \lambda \right| =1$ is necessary  is shown in proposition
\ref{prop-inv}. We sketch the construction for the case of
separable $H$. Let $\left\{ f_{n}:n=\pm 1,\pm 2,\dots \right\} $
be an orthonormal basis and  $\lambda \in \mathbb{C}$ with $\left|
\lambda \right| =1$. Define $A$ as the bilateral shift
$Af_{n}=f_{n-1}$, and $B$ as a multiplication operator via
$Bf_{n}=\lambda ^{n}f_{n}$. Both $A$ and $B$ are unitary, and $%
ABf_{n}=\lambda ^{n}f_{n-1}=\lambda BAf_{n}$. If  $\lambda $ is
not a primitive root of unity then $\sigma \left( B\right)
=\left\{ z\in \mathbb{C}:\left| z\right| =1\right\} $. Otherwise,
$B$ has finite periodic spectrum with uniform, infinite
multiplicity.

\vspace{3pt}


\noindent{\bf Example 4.} Halmos shows in fact a much stronger
result: any unitary operator in an infinite-dimensional Hilbert
space can be realized as $ABA^{-1}B^{-1}$ (Problem 239). Statement
(ii) of theorem \ref{theo1} entails, however, that the relation
$AB=UBA$ cannot be realized for every unitary $U$ if $A,B$ are
assumed to be self-adjoint: if $U=\lambda I$ then $\lambda $ can
only be $+1$ or $-1$. Theorem \ref{theo3} describes constraints on
a pair of self-adjoint operators that they have to satisfy in
order to commute up to a unitary factor.

Starting with a result of Frobenius in 1911, various authors have
studied multiplicative commutators $C=ABA^{-1}B^{-1}$ with the
property that $AC=CA$ for different types of $n\times n$ matrices
$A,B$; for a recent sample and survey, see Shapiro (1997).

\vspace{3pt}


\noindent{\bf Example 5.} Let $\left\{ f_{n}:n=0,\pm 1,\pm
2,\allowbreak\dots \right\} $ be an orthonormal basis of a
separable Hilbert space, and define $Bf_{n}=\beta _{n}f_{n}$, and
$Af_{0}=f_{-1}$, $Af_{k}=0$
for $k\neq 0$. Then $ABf_{n}=\beta _{n}\delta _{n0}f_{n-1}$, and $%
BAf_{n}=\beta _{n-1}\delta _{n0}f_{n-1}$. So $AB=\lambda BA$ is satisfied
exactly when $\beta _{0}=\lambda \beta _{-1}$. Hence every nonzero value of $%
\lambda $ can be realized by a suitable choice of $\beta
_{0}/\beta _{-1}$. In this example $A$ is nilpotent while $B$ is
normal but otherwise rather arbitrary. Note that the conditions of
lemma \ref{pro0} are violated.

Starting with the operator $B$ defined in the last paragraph, one can
investigate constraints on its spectrum and the possible operators $A$ that
give a realization of the relation $AB=\lambda BA$. In fact we obtain $%
B\left( Af_{n}\right) =\left( \beta _{n}/\lambda \right) \left(
Af_{n}\right) $, so that whenever $Af_{n}\neq 0$, this vector must
be an eigenvector of $B$ associated with some eigenvalue, say
$\beta _{m}$. It follows that $Af_{n}=\alpha _{mn}f_{m}$ and
$\beta _{m}=\beta _{n}/\lambda $ for some\ $m$. This spectral
property of $B$ must hold for all eigenvalues of $B$ whose
associated eigenvectors do not lie in the kernel of $A$, and $A$
sends the set of eigenvectors of $B$ to a subset of itself.

\vspace{3pt}


\noindent{\bf Example 6.} Let $D$ be the closed unit disk in
$\mathbb{C}$ and let $\beta :D\rightarrow D$ be a continuous
bijective map. In $H=L^{2}(D,d^{2}z)$ define $Bf(z)=\beta
(z)f(z)$, $f\in H$. This multiplication operator is normal and has
spectrum equal to $D$, and $D$ hasthe property $D=\lambda D$ for
any $\lambda =\re^{\ri\theta }$. We define $Af(z)=f\left( \beta
^{-1}\left( \lambda \beta \left( z\right) \right) \right) $, with
$A^{-1}f(z)=f\left( \beta
^{-1}\left( \lambda ^{-1}\beta \left( z\right) \right) \right) $. Then $%
ABA^{-1}f(z)=\lambda Bf\left( z\right) $. This example illustrates
proposition \ref{prop-inv}.

\vspace{3pt}


\noindent{\bf Example 7.} Next we give some simple
finite-dimensional matrix examples illustrating that any value of
$\lambda \neq 0$ may occur, with a suitable choice of $A$ for a
particular given nilpotent $B$, again in line with lemma
\ref{pro0}. Let
\begin{equation*}
B=\left(
\begin{array}{ll}
0 & 0 \\
1 & 0
\end{array}
\right) ,\;A=\left(
\begin{array}{ll}
x & 0 \\
y & \lambda x
\end{array}
\right) .
\end{equation*}
$A$ is invertible provided that $\lambda x\neq 0$. We have
\begin{equation*}
AB=\left(
\begin{array}{ll}
0 & 0 \\
\lambda x & 0
\end{array}
\right) =\lambda BA.
\end{equation*}
The same result is obtained for
\begin{equation*}
B=\left(
\begin{array}{lll}
0 & 0 & 0 \\
1 & 0 & 0 \\
0 & 1 & 0
\end{array}
\right) ,\;A=\left(
\begin{array}{ccc}
x & 0 & 0 \\
y & \lambda x & 0 \\
z & \lambda y & \lambda ^{2}x
\end{array}
\right) .
\end{equation*}
These examples are related to finite-dimensional realizations of $%
U_{q}\left( \frak{sl}\,\left( 2\right) \right) $ which are presented next.

It is a fact that on non-trivial finite-dimensional $U_{q}\left( \frak{sl}%
\,\left( 2\right) \right) $ modules, the endomorphisms associated
to $E$ and $F$ are nilpotent. On the unique, up to isomorphism,
simple $\left( n+1\right) $-dimensional $U_{q}\left(
\frak{sl}\,\left( 2\right) \right) $ module, a realization is
(Kassel 1995):
\begin{equation*}
E=\epsilon \left(
\begin{array}{ccccc}
0 & \left[ n\right] & 0 & \cdots & 0 \\
0 & 0 & \left[ n-1\right] & \cdots & 0 \\
\vdots & \ddots & \ddots & \ddots & \vdots \\
0 & 0 & \ddots & \ddots & \left[ 1\right] \\
0 & 0 & \cdots & 0 & 0
\end{array}
\right) ,\;\;F=\left(
\begin{array}{ccccc}
0 & 0 & \cdots & 0 & 0 \\
\left[ 1\right] & 0 & \cdots & 0 & 0 \\
0 & \left[ 2\right] & \ddots & 0 & 0 \\
\vdots & \ddots & \ddots & \ddots & \vdots \\
0 & 0 & \cdots & \left[ n\right] & 0
\end{array}
\right) ,
\end{equation*}
\begin{equation*}
K=\epsilon \left(
\begin{array}{ccccc}
q^{n} & 0 & \cdots & 0 & 0 \\
0 & q^{n-2} & \cdots & 0 & 0 \\
\vdots & \ddots & \ddots & \vdots & \vdots \\
0 & 0 & \cdots & q^{-n+2} & 0 \\
0 & 0 & \cdots & 0 & q^{-n}
\end{array}
\right)
\end{equation*}
where $\epsilon =\pm 1$ and $\left[ m\right] =\left( q^{m}-q^{-m}\right)
/\left( q-q^{-1}\right) =q^{m-1}+q^{m-2}+\cdots +q^{-m+2}+q^{-m+1}$. For $%
n=2 $ and $\epsilon =1$, an isomorphic three-dimensional
realization is given by Jantzen (1996):
\begin{equation*}
E=\left(
\begin{array}{ccc}
0 & q+q^{-1} & 0 \\
0 & 0 & q+q^{-1} \\
0 & 0 & 0
\end{array}
\right) ,\;\;F=\left(
\begin{array}{lll}
0 & 0 & 0 \\
1 & 0 & 0 \\
0 & 1 & 0
\end{array}
\right) ,
\end{equation*}
\begin{equation*}
K=\left(
\begin{array}{lll}
q^{2} & 0 & 0 \\
0 & 1 & 0 \\
0 & 0 & q^{-2}
\end{array}
\right)
\end{equation*}
The three-dimensional matrices $A$ and $B$ in the foregoing
example may, by the choices $y=z=0$ and $x=q^{2}$, $\lambda
=q^{-2}$, be identified with $K$ and $F$, respectively.

\section{Concluding remark}

The results of \S2 can be transferred literally into the context
of von Neumann algebras. The problem treated in this paper arises
from an operator algebraic question in the context of quantum
measurement theory (Busch \& Singh 1998). For positive bounded
operators $A,B$, define positive linear maps on the space of
bounded operators, $X\mapsto AXA$, $X\mapsto BXB$. These maps
commute, $ABXBA=BAXAB$, exactly when $AB=\lambda BA$ for some complex $%
\lambda $. If $AB=O$, then $A$ and $B$ commute. Assume $AB\neq O$.
Upon taking adjoints, $AB=\overline{\lambda }^{-1}BA$, so $\left|
\lambda \right| =1$. Theorem \ref{theo1} asserts that $\lambda
=1$. If these maps describe state transformations due to local
measurements in spacelike separated regions of Minkowski
spacetime, this result ensures that the description of state
changes is Lorentz invariant (Busch 1999). More generally, if
there is a unitary operator $U$ such that $UABXBA=UBAXABU^{\ast }$
for all $X$, then as a consequence of theorem \ref{theo3},
$AB=BA$.

\begin{acknowledgements}
We wish to thank an anonymous referee for providing an additional
angle on the subject of this paper that we were not fully aware of
in the course of writing, and for helpful organizational
suggestions.
\end{acknowledgements}

\label{lastpage}
\end{document}